\documentclass[11pt]{article} 

\usepackage[colorlinks=true,linkcolor=blue,urlcolor=blue]{hyperref}
\usepackage{dsfont}

\usepackage{amsfonts,latexsym,amstext} 
\usepackage{amsmath} 
\usepackage{amssymb} 
\usepackage{comment} 
\usepackage{amsthm} 
\usepackage[latin1]{inputenc} 
\usepackage{color}

\theoremstyle{plain} 
\newtheorem{theorem}{Theorem}[section] 
\newtheorem{lemma}[theorem]{Lemma} 
\newtheorem{proposition}[theorem]{Proposition} 
\newtheorem{corollary}[theorem]{Corollary} 
 \def\R{\rz R}

\def\P{\rz P}
\def\Q{\rz Q}

\def\shd{{\cal D}}

\def\shs{{\cal S}}

\newtheorem{definition}[theorem]{Definition} 
 
\theoremstyle{remark} 
\newtheorem{remark}[theorem]{Remark} 

\allowdisplaybreaks 
 
\makeatletter 
\@addtoreset{equation}{section} 
 
\makeatother 
 
\def\sqr#1#2{{\vcenter{\vbox{\hrule height .#2pt \hbox{\vrule 
 width .#2pt height#1pt \kern#1pt \vrule 
width .#2pt} \hrule height .#2pt}}}}

\def\ds{\begin{displaystyle}} 
\def\eds{\end{displaystyle}} 
 
\def\<{\langle } 
\def\>{\rangle }

\def\R{\mathbb R}

\def\P{\mathbb P} 
\def\Q{\mathbb Q}

\DeclareMathAlphabet{\mathonebb}{U}{bbold}{m}{n}                           %
\newcommand{\one}{\ensuremath{\mathonebb{1}}}                               

\title{
Characteristics and It\^o's formula for weak Dirichlet processes: an equivalence result}

\author{ 
Elena Bandini 
\thanks{Universit\`a di Bologna, Dipartimento di Matematica, Piazza di Porta S. Donato 5, 40126 Bologna, Italy; e-mail: elena.bandini7@unibo.it.} 
\and Francesco Russo 
\thanks{ENSTA Paris, Unité de Mathématiques Appliquées, 
  Institut Polytechnique de Paris, 
828, boulevard des Mar\'echaux, F-91120 Palaiseau, France; e-mail: francesco.russo@ensta-paris.fr.} 
} 
\date{}
 
\begin{document} 
  
\allowdisplaybreaks 
\maketitle 
 
\begin{abstract} The main objective consists in generalizing  a well-known
  It\^o formula of J. Jacod and A. Shiryaev: 
  given a càdlàg process $S$,
  there is an equivalence between the fact that $S$
  is a semimartingale with given characteristics
  $(B^k, C, \nu)$ and  a It\^o formula type expansion
  of $F(S)$, where $F$ is a bounded function of class $C^2$. This result connects weak solutions of path-dependent SDEs and related martingale problems.
  We extend this to the case when $S$ is a weak Dirichlet
  process.  
  A second aspect of the paper consists in discussing some untreated features of stochastic calculus 
  for finite quadratic variation processes. 

\end{abstract} 

{\bf Key words:}  Martingale problem, Itô formula,
weak Dirichlet process.
 \\
 {\small\textbf{MSC 2020:}  
60H48; 60H05; 60H10.
}

\section{Introduction}


The motivating formula is the one of J. Jacod and
A. Shiryaev concerning a generic càdlàg process $X$.
Let $k$ be a truncation function cutting large jumps, i.e., a bounded function such that $k(x)=x$ in a neighborhood of zero. 
This says that $X$ is a semimartingale with
characteristics $(B^k, C, \nu)$ if and only if
for every $F \in C^{1,2}_b([0,T] \times \R)$, 
\begin{align} \label{JacodEquiv} 
  &F(\cdot, X_{\cdot}) - F(0,X_0) -\int_0^\cdot \partial_s F(s, X_s) ds- \frac{1}{2} \int_0^{\cdot}  \partial_{xx}^2 F(s,X_s) \,d(C\circ X)_s-
    \int_0^{\cdot}  \partial_{x} F(s,X_{s-}) \,d (B^{k}\circ X)_s\notag\\
&-  \int_{]0,\,\cdot]\times \R} (F(s,X_{s-} + x) -F(s,X_{s-})-k(x)\, \partial_{x} F(s,X_{s-}))\,(\nu\circ X)(ds\,dx)
\end{align} 
is a local martingale, i.e. a It\^o formula expansion.
This equivalence was stated in Theorem  2.42, Chapter II, in \cite{JacodBook}, where $F$ is time homogeneous, but can be easily extended to the non-homogeneous case, see also Appendix \ref{App:A}. 
 Decomposition \eqref{JacodEquiv} can 
be seen as a martingale problem formulation
for a class of function $F \in C^{1,2}_b([0,T] \times \R)$ for which
$F(\cdot, X)$ is a prescribed special semimartingale. 


Given a filtration $(\mathcal F_t)$, a  (càdlàg) $(\mathcal F_t)$-weak Dirichlet process is defined as the sum of a local
$(\mathcal F_t)$-martingale $M$ and an $(\mathcal F_t)$-{\it martingale orthogonal} process $A$
which means that $[A,N] = 0$, where $N$
is a generic continuous $(\mathcal F_t)$-martingale. Typical examples of martingale orthogonal processes are by definition  purely discontinuous martingales and bounded variation processes, see Proposition 2.14 in \cite{BandiniRusso1}. 
In particular, if $X$ is an $(\mathcal F_t)$-semimartingale, then $X$ is an $(\mathcal F_t)$-weak Dirichlet process.
If one forces $M$ to be continuous,  
then the decomposition $X = M + A$
(fixing $A_0= 0$) is unique, see Proposition 3.2 in \cite{BandiniRusso_RevisedWeakDir}.
In this case the continuous local martingale component $M$ is denoted
by $X^c$.
If $A$ is predictable then it is called $(\mathcal F_t)$-{\it special weak Dirichlet process}. The filtration $(\mathcal F_t)$ will be often omitted when it is self-explanatory or it is the canonical filtration $(\mathcal F^X_t)$ associated with   the process $X$. 

The concept of continuous weak Dirichlet process was introduced in \cite{er2}, and 
further analyzed and extended to the multidimensional case in \cite{gr}; a first application to stochastic
control was performed in \cite{gr1}.
In particular therein the authors investigated
the stability of weak Dirichlet processes
through $C^{0,1}$-transformations.
Indeed, if $X$ is an $\R^d$-valued weak Dirichlet  process admitting all its mutual covariations, and $F \in C^{0,1}([0,T] \times \R^d)$, 
then $F(\cdot, X)$ is still a weak Dirichlet process.
In the jump case, \cite{cjms} introduced the notion of special Dirichlet process
generalizing the notion of special semimartingale.
In that framework the general notion of weak Dirichlet process appeared
in \cite{BandiniRusso1} and was applied to the BSDEs theory in \cite{BandiniRusso2}.
The aforementioned stability property has been extended in
\cite{BandiniRusso_RevisedWeakDir}
to the case when $X$ is a jump process.
A generalization of the  $C^{0,1}$-stability 
in the continuous framework, but $u$ being a path-dependent functional, 
has been performed in \cite{BLT} with application
to finance. A survey on weak Dirichlet processes is provided in Chapter 15
of \cite{Russo_Vallois_Book}.

The notion of characteristics, well known for semimartingales,  has been extended in 
\cite{BandiniRusso_RevisedWeakDir}
to a generic weak Dirichlet process as follows. Cutting large jumps via a truncation function $k$,
the process
\begin{equation}\label{specXI}
X_t - \sum_{s \le t} (\Delta X_s - k(\Delta X_s)),\quad  t \geq 0,
\end{equation}
is a special weak Dirichlet process. If $(\mathcal F_t)=(\mathcal F^X_t)$, 
then the process \eqref{specXI} can be decomposed as
$$
B^k \circ X + X^c + k(x) \star (\mu^X - \nu\circ X),
$$  with 
$[X^c,X^c] = C \circ X$ and $\nu$  a random measure
with $\nu \circ X$ the compensator of the jump counting
measure $\mu^X$. Therefore, $X$ fulfills 
the equation
\begin{equation}\label{decIntr}
X = X^c + k(x) \star (\mu^X- \nu\circ X)+ B^k \circ X + (x - k(x)) \star \mu^X.
\end{equation}

Our main Theorem
\ref{T:to_prove} provides a generalization of the  Jacod-Shiryaev equivalence
theorem (see \eqref{JacodEquiv}) when $X$ is not necessarily a semimartingale but a weak Dirichlet process, see Definition \ref{finitequadvar}. Namely it states that a finite quadratic variation process   
X is a weak Dirichlet process with local characteristics $(B^k, C, \nu)$ if and only if, for each bounded function $F$ of class $C^{1,2}([0,T] \times \R)$, 
$
\int_{]0,\cdot]}  \partial_x F(s,X_{s}) \,d^- (B^{k}\circ X)_s
$
is a martingale orthogonal process, 
 and the process 
\begin{align}\label{f2}
&F(\cdot, X_{\cdot}) - F(0, X_0) -\int_0^\cdot \partial_s F(s, X_s) ds- \frac{1}{2} \int_0^{\cdot} \partial^2_{xx} F(s,X_s) \,(d(C \circ X)_s +d[B^k \circ X, B^k \circ X]^c_s)\notag\\
&-  \int_{]0,\,\cdot]} \partial_x F(s,X_{s}) \,d^{-} (B^{k}\circ X)_s\notag\\
&-  \int_{]0,\,\cdot]\times \R}(F(s,X_{s-} + x) -F(s,X_{s-})-k(x)\,\partial_x F(s,X_{s-}))\,(\nu\circ X)(ds\,dx)
\end{align} 
is an $({\mathcal F_t^X})$-local martingale, where $ \int_{]0,\,\cdot]}  \partial_x F(s,X_{s}) \,d^{-} (B^{k}\circ X)_s$ is the forward integral in the sense of \cite{rv95} (see Definition \ref{D_ucp_integral}), see  \cite{BandiniRusso1}, and also \cite{nunno-book}.  
For proving \eqref{f2}, a fundamental tool  of independent interest is the following (see Theorem \ref{L:Bfinitevar}): if  $X$ is an $(\mathcal F^X_t)$-weak Dirichlet process with finite quadratic  variation with decomposition \eqref{decIntr}, 
	then $B^{k}\circ X$ is a finite quadratic variation process, with  
  \begin{align}\label{brack_B2_Intr}
 	[B^{k}\circ X, B^{k}\circ X]_t 
 	&=[X,X]^c_t -[X^c,X^c]_t 
 	+\sum_{s \leq t} \Big|\int_\R k(x)  (\nu\circ X)(\{s\} \times dx)\Big|^2.  
 \end{align} 
This  in particular extends the results when $X$ is  semimartingale or a Dirichlet process: in those cases, $[B^{k}\circ X, B^{k}\circ X]_t  = \sum_{s \leq t}|\Delta (B^{k}\circ X)|^2= \sum_{s \leq t} |\int_\R k(x)  (\nu\circ X)(\{s\} \times dx)|^2$ and therefore $[X,X]^c_t =[X^c,X^c]_t$, see Remark \ref{R:3.5}.

An important  aspect of Theorem
\ref{T:to_prove} is that it prolongates the classical equivalence between weak solutions of stochastic differential equations and martingale problems in the framework of continuous Markov
processes, see \cite{Stroock_Varadhan}.

 Formulation \eqref{f2} fits the one of Theorem 4.3 in \cite{BandiniRusso_RevisedWeakDir}, where we specify $F(\cdot, X)$ for $F$ belonging to some domain $\shd_\shs$ which
is a subspace of $C^{0,1}([0,T] \times \R)$.  In the present  case,
$\shd_\shs$  is the space of bounded functions of class $C^{1,2}([0,T] \times   \R)$, and  the processes $F(\cdot, X)$  are no longer special semimartingales but special weak Dirichlet processes. 
In \cite{BandiniRusso_RevisedWeakDir} we were specially interested in domains for  which the processes  $F(\cdot, X)$ were still semimartingales, see Definition 4.12 in \cite{BandiniRusso_RevisedWeakDir}; this included the framework of $X$ being a solution of a stochastic differential equation with jump
 and singular (distributional) drift, see Theorem 4.1 in \cite{BandiniRusso_DistrDrift}.

A second aspect of the paper consists in  discussing some properties of the covariation processes. 
We give an explicit expression for the covariation of   two c\`adl\`ag  finite quadratic variation processes, provided that  the continuous component of the quadratic variation of one of them  vanishes, see Lemma \ref{L:YZ}. 
We also provide a stability result  for $C^1$ transformations of finite quadratic variation processes, see Lemma  \ref{L:app1}.



 

The paper is organized as follows. Section \ref{s:ESC} is devoted to some useful results concerning stochastic calculus via regularization for jump processes,
including the related properties to the covariance.  In Section
\ref{S:WDP} we formulate some basic recalls on weak Dirichlet processes.
Section \ref{chp4} is devoted to the statement and the proof of the main Theorem \ref{T:to_prove}.

\section{Elements of stochastic calculus via regularization for jump processes}
\label{s:ESC}

\subsection{Preliminaries}
In the whole article, we are given a fixed maturity $T>0$  and a  probability space $(\Omega,\mathcal{F},\P)$.

We will consider the space of functions
$ 
u: [0,T] \times \R \rightarrow \R$, $(t,x)\mapsto u(t,x)$, 
 which are of class  $C^{0,1}$ or  $C^{1,2}$.
 $C^{0,1}_b$  (resp. $C^{1,2}_b$) stands for the class of bounded
 functions which belong to  $C^{0,1}$ (resp. $C^{1,2}$).
 $C^{0,1}$ is equipped with the topology of uniform convergence on each compact
 of $u$ and $\partial_x u$. 
$C^0$ (resp. $C^0_b$) will denote the space of continuous functions (resp. continuous and bounded  functions) on $\R$ equipped with the topology of uniform convergence on each compact (resp. equipped with the topology of uniform convergence), while   $C^1$ (resp. $C^2$) will be the space of continuously differentiable (twice continuously differentiable) functions $u:\R\rightarrow \R$.  $C^{1}_b$  (resp. $C^{2}_b$) stands for the class of bounded
 functions which belong to  $C^{1}$ (resp. $C^{2}$). 


The concept of random measure   
will be extensively used 	throughout   the paper. 
For a  detailed  discussion on this topic  and the unexplained  notations,
we refer to
Chapter I and Chapter II, Section 1, in \cite{JacodBook}, Chapter III in \cite{jacod_book},  and  Chapter XI, Section 1, in \cite{chineseBook}.
In particular, if $\mu$ is a random measure on $[0,\,T]\times \R$, for any measurable real function $H$ defined on $\Omega \times [0,\,T] \times \R $, one denotes 
$$
H \star \mu_t:= \int_{]0,\,t] \times \R} H(\cdot, s,x) \,\mu(\cdot, ds \,dx),
$$
when the stochastic integral in the right-hand side is defined
(with possible infinite values).


 Let $X$ be an adapted c\`adl\`ag process. We denote by $\Delta X$, with  $\Delta X_t = X_t - X_{t-}$, the corresponding jump process. 
 We set the 
 corresponding  jump measure  $\mu^X$ by
\begin{equation}\label{jumpmeasure}
\mu^X(dt\,dx)= \sum_{s >0} \one_{\{\Delta X_s \neq 0\}}\, \delta_{(s, \Delta X_s)}(dt\,dx). 
\end{equation}
In this case, $H \star \mu^X_t= \sum_{0 < s \leq t} H(\cdot, s,\Delta X_s)$.
We denote by  $\nu^X= \nu^{X, \P}$ the compensator of $\mu^X$,
         see \cite{JacodBook}, Theorem 1.8, Chapter II.
         The dependence on $\P$ will be omitted when self-explanatory.


\begin{definition}\label{D_ucp_integral} 
Let $X$ be a c\`adl\`ag process and $Y$ be a  process belonging to $L^1([0,\,T])$ a.s. 

 $\int_{]0,\,\cdot]}\,Y_s \, d^{-} X_s$ denotes the forward integral of $Y$ with respect to $X$, i.e., the u.c.p. limit, whenever it exists, of 
 \begin{eqnarray*}
	\int_0^t Y(s)\,\frac{X((s+\varepsilon)\wedge t)-X(s)}{\varepsilon}\,ds, \quad t \in [0,T]. 
\end{eqnarray*}
\end{definition} 
\begin{remark}\label{R:12}
	Let $(Y, Z)$, $(Y',Z')$ be two pair of processes and $\Omega_0 \subset \Omega$ be an event such that 
	$$
	Y_t \one_{\Omega_0}= Y'_t \one_{\Omega_0},\quad Z_t \one_{\Omega_0}= Z'_t \one_{\Omega_0}, \quad \forall t \in [0,T], 
	$$
	where the equality holds in the sense of indistinguishability. 
	Then
	$$
	\int_{]0,\cdot]} Y_s d^- Z_s= \int_{]0,\cdot]} Y'_s d^- Z'_s
	$$
	in the sense that, if an integral exists, then the other exists and they are equal. 
\end{remark}
\begin{definition}\label{finitequadvar}
For two c\`adl\`ag processes $X$ and $Y$,  we define the covariation of $X$ and $Y$, denoted $[X, Y]$, as the   u.c.p. limit (if it exists) of  
\begin{equation} 
\label{Appr_cov_ucpI}	 
[X,Y]^{\varepsilon}(t):= \,\int_0^t\,\frac{(X((s+\varepsilon)\wedge t)-X(s))(Y((s+\varepsilon)\wedge t)-Y(s))}{\varepsilon}\,ds, \quad t \in [0,T]. 
\end{equation} 
A c\`adl\`ag process $X$ will be called  a finite quadratic variation process 
 whenever  $[X,X]$     
exists. 
\end{definition} 
By Lemma 2.10 in \cite{BandiniRusso1}, we know that
\begin{equation} \label{QVC}
[X,X] = [X, X]^c  + \sum_{s \leq \cdot} |\Delta X_s|^2,
\end{equation}
where $[X,X]^c$ is the continuous component of $[X,X]$.

\begin{remark}\label{R:1.1} 
By Proposition 1.1 in \cite{rv95}, 
if  $X, Y$ are two c\`adl\`ag semimartingales
and $H$ is a  c\`adl\`ag adapted  process we have the following.  
\begin{itemize} 
\item[(i)] $[X,Y]$ exists and it is the usual bracket. 
\item[(ii)] $\int_{]0,\,\cdot]} H\, d^{-}X$  is the usual stochastic integral 
 $\int_0^\cdot H_{s-} dX_{s}$. 
\end{itemize} 
\end{remark}

\subsection{New  technical results}

\begin{proposition}\label{P:3.1}
Let $\varphi: [0,T] \times \R \times \Omega \rightarrow \R$ be a measurable function such that 
$$
|\varphi(s, x)| \star \mu^X_T =\sum_{s \leq T} |\varphi(s, \Delta X_s)| <+ \infty \quad \textup{a.s.}
$$
Set  $$
\tilde A_t: = \varphi(s, x) \star \mu^X_t.  
$$
Then, for every  c\`adl\`ag process $H$, 
$$
\int_{]0,\cdot]}H_s d^{-} \tilde A_s= H_{s-}  \,\varphi(s, x)\star \mu^X. 
$$

\end{proposition}
\proof 
By Remark \ref{R:1.1}-(ii), 
$$
\int_{]0,\cdot]}H_s d^{-} \tilde A_s= \int_0^\cdot H_{s-} \, d\tilde A_s. 
$$  
Indeed,  the right-hand side  is a Lebesgue integral, so no condition of adaptability   on $H$ is required. On the other hand, 
$$
\int_0^\cdot H_{s-} \, d\tilde A_s = \sum_{s \leq \cdot}H_{s-} \, \Delta\tilde A_s= \sum_{s \leq \cdot}H_{s-} \, \varphi(s, \Delta X_s)=H_{s-} \,\varphi(s,x) \star \mu^X. 
$$
\endproof

\begin{lemma}\label{L:YZ}
	Let $Y, Z$ be two c\`adl\`ag finite quadratic variation processes. Suppose that  $[Y,Y]^c=0$. Then 
	$$
	[Y,Z]= \sum_{t \leq \cdot}\Delta Y_t \Delta Z_t.
	$$
\end{lemma}
\proof
Given a sequence $(\varepsilon_n)$ converging to zero, we need to show the existence of $(n_k)$ so that 
\begin{align*}
	[Y,Z]^{\varepsilon_{n_k}}&\underset{k \rightarrow + \infty}{\longrightarrow} \sum_{t \leq \cdot}\Delta Y_t \,\Delta Z_t,\quad \textup{a.s. uniformly.}
\end{align*}
By extraction of subsequences, we can suppose that for almost all $\omega$,  there is a subsequence of $(\varepsilon_n)$, still denoted by the same letter,   such that 
\begin{align*}
	[Y,Y]^{\varepsilon_n}&\underset{n \rightarrow + \infty}{\longrightarrow}[Y,Y]= \sum_{t \leq \cdot}|\Delta Y_t|^2, \quad \textup{uniformly},\\
[Z,Z]^{\varepsilon_n}&\underset{n \rightarrow + \infty}{\longrightarrow} [Z,Z],\quad \textup{uniformly}.
\end{align*}
Let us thus fix a realization of $\omega \in \Omega$. Let $(t_i)_{i \geq 1}$ be the sequence including the jumps of $Y(\omega)$ and $Z(\omega)$, obviously in $]0,T]$.  In the sequel we will omit  the dependence on $\omega$. 

 By \eqref{QVC}, we can take 
 $\gamma >0$ and $N= N(\gamma)$ such that 
\begin{align*}
	\sum_{i=N+1}^\infty |\Delta Z_{t_i}|^2+ \sum_{i=N+1}^\infty |\Delta Y_{t_i}|^2 \leq \gamma. 
\end{align*}
We proceed similarly as for the proof of Proposition 2.14 in \cite{BandiniRusso1}. By renumerating  increasingly the set $(t_i)_{i=1}^N$ and setting $t_0=0$, we define 
\begin{align}
	A(\varepsilon, N) &= \bigcup_{i=1}^N \,]t_{i}-\varepsilon, t_i], \label{Aset}\\
		B(\varepsilon, N) &= \bigcup_{i=1}^N \,]t_{i-1}, t_i- \varepsilon]= [0,T]\setminus A(\varepsilon, N),\label{Bset}
\end{align}
with $\varepsilon < \inf_{i=1,..., N}|t_i- t_{i-1}|$. We decompose 
\begin{equation}\label{3bis}
\frac{1}{\varepsilon} \int_0^s(Y_{(t+ \varepsilon)\wedge s}-Y_t)(Z_{(t+ \varepsilon)\wedge s}-Z_t)dt - \sum_{t \leq s} \Delta Y_t \Delta Z_t
\end{equation}
into
\begin{equation}\label{3}
I_A^{Y, Z}(\varepsilon, N, s)+I_{B_1}^{Y, Z}(\varepsilon, N, s)+I_{B_2}^{Y, Z}(N, s), 
\end{equation}
where 
\begin{align*}
	I_{A}^{Y, Z}(\varepsilon, N, s)&=\frac{1}{\varepsilon} \int_{]0,s] \cap A(\varepsilon, N)}(Y_{(t+ \varepsilon)\wedge s}-Y_t)(Z_{(t+ \varepsilon)\wedge s}-Z_t)dt- \sum_{i=1}^N\one_{]0,s]}(t_i) \Delta Y_{t_i} \Delta Z_{t_i},\\
	I_{B1}^{Y, Z}(\varepsilon, N, s)&=\frac{1}{\varepsilon} \int_{]0,s] \cap B(\varepsilon, N)}(Y_{(t+ \varepsilon)\wedge s}-Y_t)(Z_{(t+ \varepsilon)\wedge s}-Z_t)dt, \\
I_{B2}^{Y, Z}(N, s)&=- \sum_{i=N+1}^\infty\one_{]0,s]}(t_i) \Delta Y_{t_i} \Delta Z_{t_i}.
\end{align*}
In order to prove that 
\begin{equation}\label{conv_A}
		I_A^{Y, Z}(\varepsilon_n, N, \cdot)\underset{n \rightarrow + \infty}{\rightarrow} 0\quad \textup{uniformly},
\end{equation}
by bilinearity it is enough to show that $I_A^{\eta, \eta}(\varepsilon_n, N, \cdot )$ goes to zero, uniformly,  for $\eta= Y, Z, Y+Z$.
Now, by Lemma 2.11 of \cite{BandiniRusso1},
\begin{align*}
I_A^{\eta, \eta}(\varepsilon_n, N, s)
+ \sum_{i=1}^N\one_{]0,s]}(t_i) (\Delta \eta_{t_i})^2&=  
		\sum_{i=1}^N \frac{1}{\varepsilon} \int_{t_i - \varepsilon}^{t_i} \one_{]0,\,s]}(t)\,\phi(\eta_{(t+ \varepsilon)\wedge s},\eta_t)\,dt \\
		&\overset{\varepsilon \rightarrow 0}{\longrightarrow}\,\, \sum_{i=1}^N  \one_{]0,\,s]}(t_i)\, \phi(\eta_{t_i},\eta_{t_{i}-}),
	\quad \textup{uniformly in }\,\, s \in [0,\,T], 
\end{align*}
with $\phi(x_1, x_2)= (x_1 - x_2)^2$. 
This implies that
$
I_A^{\eta, \eta}(\varepsilon_n, N, \cdot)$ goes  to zero, uniformly.

Since  $[Y,Y]^c=0$, taking into account \eqref{QVC} we have \begin{equation}\label{conv1}
\frac{1}{\varepsilon} \int_0^s(Y_{(t + \varepsilon_n)\wedge s}-Y_t)^2 dt - \sum_{t \leq s}|\Delta Y_t|^2 \underset{n \rightarrow + \infty}{\longrightarrow} 0, \quad \textup{uniformly in }\,\, s \in [0,\,T]. 
\end{equation}
On the other hand, 
\begin{equation}\label{B1}
|I_{B2}^{Y, Y}(N, s)|\leq \sum_{i=N+1}^\infty  |\Delta Y_{t_i} |^2 \leq \gamma. 
\end{equation}
Collecting  \eqref{conv1},  \eqref {conv_A} with $Z=Y$ and  \eqref{B1}, it follows from \eqref{3bis} and \eqref{3} both for $Z=Y$ that 
\begin{equation}\label{IB1_limsup}
\limsup_{n \rightarrow +\infty}
|I_{B1}^{Y, Y}(\varepsilon_n, N, s)| \leq \gamma. 
\end{equation}

We come back to the estimate of \eqref{3bis}. The absolute value of \eqref{3bis},  with $\varepsilon = \varepsilon_n$, is bounded by 
\begin{align}\label{5}
	&|I_A^{Y, Z}(\varepsilon_n, N, s)|+|I_{B_1}^{Y, Z}(\varepsilon_n, N, s)|+|I_{B_2}^{Y, Z}(N, s)| \notag\\
		&\leq |I_A^{Y, Z}(\varepsilon_n, N, s)| + \sqrt{I_{B1}^{Y, Y}(\varepsilon_n, N, s)} \sqrt{I_{B1}^{Z,Z}(\varepsilon_n, N, s)} 
	+\sqrt{\sum_{i=N+1}^\infty  |\Delta Y_{t_i} |^2}\sqrt{\sum_{i=N+1}^\infty  |\Delta Z_{t_i} |^2}\notag\\
	&\leq |I_A^{Y, Z}(\varepsilon_n, N, s)|+ \sqrt{I_{B1}^{Y, Y}(\varepsilon_n, N, s)} \, \sqrt {[Z,Z]^{\varepsilon_n}_T} + \gamma. 
\end{align} 
Taking the $\limsup_{n \rightarrow \infty}$ in   \eqref{5}, taking into account \eqref{conv_A} and \eqref{IB1_limsup},  we get  
\begin{align*}
	&\limsup_{n \rightarrow \infty}\, (|I_A^{Y, Z}(\varepsilon_n, N, s)|+|I_{B_1}^{Y, Z}(\varepsilon_n, N, s)|+|I_{B_2}^{Y, Z}(N, s)|) \leq \sqrt \gamma \, \sqrt{[Z,Z]_T }+ \gamma. 
\end{align*}
Since $\gamma$ is arbitrarily chosen, we have shown that  \eqref{3bis} converges uniformly to zero. This implies that \eqref{3bis} converges u.c.p. to zero. 
\endproof


The following result of stability of finite quadratic variation processes was well understood in the context of F\"ollmer's discretizations,  but was never established in the regularization framework. 
\begin{lemma}\label{L:app1}
	\begin{enumerate}
	\item
	Let $Y= \varphi(X)$, where $\varphi: \R \rightarrow \R$ is a $C^1$ function and $X$ is a  c\`adl\`ag process of finite quadratic variation. Then
	\begin{equation*}
		[Y,Y]_t = \int_0^t (\varphi'(X_{s-})^2  d[X, X]_s^c + \sum_{s \leq t}(\Delta \varphi(X_s))^2, \quad t \in [0,T].
	\end{equation*} In particular,  $Y$ is also a  finite quadratic variation process.
	\item 	Let $Y^1= \varphi(X^1)$ and $Y^2= \phi(X^2)$, where $\varphi$ and $\phi$ are $C^1$ functions and $X^1, X^2$ are   c\`adl\`ag processes such that $(X^1, X^2)$ has all its mutual covariations. Then 
	\begin{equation*}
		[Y^1,Y^2]_t = \int_0^t \varphi'(X^1_s) \phi'(X^2_{s-}) d[X^1, X^2]_s^c + \sum_{s \leq t}\Delta \varphi(X^1_s)\,\Delta \phi(X^2_s),\quad  t \in [0,T].
	\end{equation*}
	\end{enumerate}

\end{lemma}
\proof
1. Let $t \in [0,\,T]$ and $\varepsilon \in [0,\,1]$. We expand, for $s \in [0,\,T]$,
\begin{align*}
	\varphi(X_{(s + \varepsilon)\wedge t})-\varphi(X_{s \wedge t})= I_1^\varphi(s, t, \varepsilon) (X_{(s + \varepsilon)\wedge t}-X_{s \wedge t}), 
\end{align*}
where
 \begin{align*}
   I_1^\varphi(s, t, \varepsilon)= \int_0^1 \varphi'(X_{s \wedge t}+
   a (X_{(s + \varepsilon)\wedge t}-X_{s\wedge t}))\,da.
  \end{align*}
Consequently, 
\begin{align*}
	\frac{1}{\varepsilon}(\varphi(X_{(s + \varepsilon)\wedge t})-\varphi(X_{s \wedge t}))^2 &= \frac{1}{\varepsilon}((I_1^\varphi(s, t, \varepsilon))^2-(\varphi'(X_s))^2) (X_{(s + \varepsilon)\wedge t}-X_{s \wedge t})^2 \\
	&+\frac{1}{\varepsilon}(\varphi'(X_s))^2 (X_{(s + \varepsilon)\wedge t}-X_{s \wedge t})^2. 
\end{align*}
Integrating from $0$ to $t$, we get 
\begin{align}\label{J1+J2}
	\frac{1}{\varepsilon}\int_0^t(\varphi(X_{(s + \varepsilon)\wedge t})-\varphi(X_{s}))^2 ds &= \frac{1}{\varepsilon}\int_0^t((I_1^\varphi(s, t, \varepsilon))^2-(\varphi'(X_s)^2) (X_{(s + \varepsilon)\wedge t}-X_{s})^2 ds\notag\\
	&+\frac{1}{\varepsilon}\int_0^t (\varphi'(X_s))^2 (X_{(s + \varepsilon)\wedge t}-X_{s })^2 ds\notag\\
	&=: J_1(t, \varepsilon) + J_2(t, \varepsilon).
\end{align}
We notice that, without restriction of generality, passing to a suitable subsequence, we can suppose (with abuse of notation) that
\begin{equation}\label{convbrac}
[X,X]^\varepsilon:=\frac{1}{\varepsilon} \int_0^\cdot \  (X_{(s + \varepsilon)\wedge \cdot}-X_{s})^2\, ds\underset{\varepsilon \rightarrow 0}{\rightarrow}  [X,X], \quad \textup{uniformly a.s.}
\end{equation}
Since $X$ is a finite quadratic variation process, by Lemma A.5 in \cite{BandiniRusso1}, taking into account Definition A.2 and Corollary A.4-2. in \cite{BandiniRusso1}, if $g$ is a c\`adl\`ag process then 
\begin{align}\label{B9}
	\frac{1}{\varepsilon}\int_0^t g_{s} (X_{(s + \varepsilon) \wedge t} -X_{s})^2  ds=\frac{1}{\varepsilon}\int_0^t g_{s-} (X_{(s + \varepsilon) \wedge t} -X_{s})^2  ds \underset{\varepsilon \rightarrow 0}{\rightarrow} \int_0^t g_{s-} d [X,X]_s, \quad \textup{u.c.p.}
\end{align}
Therefore, taking $g_s= (\varphi'(X_s))^2$ in  \eqref{B9}, we get 
\begin{align}\label{J2conv}
J_2(\cdot, \varepsilon) \underset{\varepsilon \rightarrow 0}{\rightarrow} \int_0^\cdot (\varphi'(X_{s-}))^2 d [X,X]_s, \quad  \textup{u.c.p.}	
\end{align}

Next step consists in proving that
\begin{align}\label{toproveJ1}
	J_1(\cdot, \varepsilon) \underset{\varepsilon \rightarrow 0}{\rightarrow} \sum_{s \leq \cdot}\Big[\Big(\int_0^1 \varphi'(X_{s-} + a \Delta X_s) da\Big)^2 - (\varphi'(X_{s-}))^2\Big](\Delta X_s)^2, \quad \textup{u.c.p.}
\end{align}
We fix a realization $\omega \in \Omega$. Proceeding as in the proof of Proposition 2.14 in \cite{BandiniRusso1}, let  $(t_i)$ be an enumeration of all the jumps of $X(\omega)$ in $[0, \, T]$. We have 
$ \sum_i (\Delta X_{t_i}(\omega))^2 < \infty$. 
Let $\gamma >0$ and $N= N(\gamma)$ such that 
\begin{align}\label{estXsquare}
	\sum_{i=N+1}^\infty  (\Delta X_{t_i}(\omega))^2 \leq \gamma^2. 
\end{align}
We decompose 
\begin{align}\label{J1}
	J_1(t, \varepsilon) &=\frac{1}{\varepsilon} \int_0^t \one_{A(\varepsilon, N)}(s) J_{1 0}(s, t,  \varepsilon) ds +\int_0^t \one_{B(\varepsilon, N)}(s) J_{1 0}(s, t,  \varepsilon) ds\notag\\
	&=:J_{1A}(t, \varepsilon, N) + J_{1B}(t, \varepsilon, N),
\end{align}
where we have denoted 
$$
J_{1 0}(s, t,  \varepsilon):=(X_{(s + \varepsilon)\wedge t}-X_{s})^2((I_1^\varphi(s, t, \varepsilon))^2-(\varphi'(X_s)^2),  
$$
and the sets $A(\varepsilon, N)$ and $B(\varepsilon, N)$ are the ones  introduced in \eqref{Aset}-\eqref{Bset}. 
By Lemma 2.11 in \cite{BandiniRusso1}, it follows that,  uniformly in $t \in [0,\,T]$,
\begin{align}\label{6}
	J_{1A}(t, \varepsilon, N) 
	&\underset{\varepsilon \rightarrow 0}{\rightarrow} \sum_{i=1}^N \one_{]0,\,t]}(t_i)(\Delta X_{t_i})^2\Big(\Big(\int_0^1 \varphi'(X_{t_i -} + a \Delta X_{t_i}) da\Big)^2- (\varphi'(X_{t_i -}))^2\Big).
	\end{align}
On the other hand,
\begin{align*}
	J_{1B}(t, \varepsilon, N) &= \sum_{i=1}^{N}\frac{1}{\varepsilon} \int_0^t \ (X_{(s + \varepsilon)\wedge t}-X_{s})^2 I^{\varphi,i}_{1 B}(s, t,  \varepsilon) \,ds,
\end{align*}  
where 
\begin{align*}
	I^{\varphi,i}_{1B}(s, t,  \varepsilon) &= \one_{]t_{i-1}, t_{i}- \varepsilon[}(s)  \Big[\Big(\int_0^1 \varphi'(X_{s \wedge t}+a(X_{(s + \varepsilon)\wedge t}-X_{s\wedge t}))\,da\Big)^2-(\varphi'(X_s))^2\Big]\\
	& =\one_{]t_{i-1}, t_{i}- \varepsilon[}(s)  \Big[\int_0^1 \varphi'(X_{s \wedge t}+a(X_{(s + \varepsilon)\wedge t}-X_{s\wedge t}))\,da-\varphi'(X_s)\Big]\cdot\\
	&\,\,\,\,\cdot\Big[\int_0^1 \varphi'(X_{s \wedge t}+a(X_{(s + \varepsilon)\wedge t}-X_{s\wedge t}))\,da+\varphi'(X_s)\Big]. 
\end{align*}
For every $i =1,..., N$, we have 
\begin{align*}
	|I^{\varphi,i}_{1 B}(s, t,  \varepsilon)| \leq 2 \sup_{y \in [X_s, X_{s +\varepsilon}]}|\varphi'(y)|\,\delta\Big(\varphi, \, \sup_i \sup_{\underset{|p-q|\leq \varepsilon}{p, q \in ]t_{i-1}, t_i[}}|X_p - X_q|\Big). 
\end{align*}
We notice that there is $\varepsilon_0$ such that, if $\varepsilon < \varepsilon_0$, $\sup_{\underset{|p-q|\leq \varepsilon}{p, q \in ]t_{i-1}, t_i[}}|X_p - X_q|\leq 3 \gamma$, 
where we have applied  
Lemma 2.12 in \cite{BandiniRusso1} to the prolongation by continuity of  $X$  to the extremities restricted to  $]t_{i-1},t_i[$. 
Therefore, for $\varepsilon < \varepsilon_0$, 
 \begin{align*}
 	|I^{\varphi,i}_{1 B}(s, t,  \varepsilon)| \leq  2 \sup_{y \in [-||X||_\infty, ||X||_\infty]}|\varphi'(y)|\,\delta(\varphi, \,3 \gamma),
 	\end{align*}
and consequently,   
\begin{align}\label{9}
	\sup_{t \in [0,\,T]}|J_{1B}(t, \varepsilon, N)| &\leq    2 \delta(\varphi, \,3 \gamma)\, \sup_{t \in [0,\,T]}[X,X]^\varepsilon_t\,\sup_{y \in [-||X||_\infty, ||X||_\infty]}|\varphi'(y)|,
\end{align}
where the latter supremum is finite by  \eqref{convbrac}.
Going back to \eqref{J1} we get  
\begin{align}\label{J1est}
	&\sup_{t \in [0,\,T]}\Big|J_1(t, \varepsilon)-\sum_{i=1}^\infty \one_{]0,\,t]}(t_i)(\Delta X_{t_i})^2\Big[\Big(\int_0^1 \varphi'(X_{t_i-} + a \Delta X_{t_i}) da\Big)^2 - (\varphi'(X_{t_i-}))^2\Big]\Big| \notag\\
	&\leq \sup_{t \in [0,\,T]}\Big|J_{1A}(t, \varepsilon, N)-\sum_{i=1}^N \one_{]0,\,t]}(t_i)(\Delta X_{t_i})^2\Big[\Big(\int_0^1 \varphi'(X_{t_i-} + a \Delta X_{t_i}) da\Big)^2 - (\varphi'(X_{t_i-}))^2\Big]\Big| \notag\\
	&+ \sum_{i=N+1}^\infty \one_{]0,\,T]}(t_i)(\Delta X_{t_i})^2\Big[\Big(\int_0^1 \varphi'(X_{t_i-} + a \Delta X_{t_i}) da\Big)^2 - (\varphi'(X_{t_i-}))^2\Big] \notag\\
	& + \sup_{t \in [0,\,T]}|J_{1B}(t, \varepsilon, N)|.
\end{align}
Taking the $\limsup_{\varepsilon \rightarrow 0}$ in \eqref{J1est}, collecting  \eqref{6} and \eqref{9},  we get 
\begin{align*}	
&\limsup_{\varepsilon \rightarrow 0}\sup_{t \in [0,\,T]}\Big|J_1(t, \varepsilon)-\sum_{i=1}^\infty \one_{]0,\,t]}(t_i)(\Delta X_{t_i})^2\Big[\Big(\int_0^1 \varphi'(X_{t_i-} + a \Delta X_{t_i}) da\Big)^2 - (\varphi'(X_{t_i-}))^2\Big]\Big| \notag\\
	&\leq 2\sum_{i=N+1}^\infty \one_{]0,\,T]}(t_i)(\Delta X_{t_i})^2 \sup_{y \in [-||X||_\infty, ||X||_\infty]}|\varphi'(y)|^2
	\\&+ 2 \sup_{\varepsilon < \varepsilon_0}  \sup_{t \in [0,\,T]}[X,X]_t^\varepsilon\,\delta(\varphi, \,3 \gamma)\, \sup_{y \in [-||X||_\infty, ||X||_\infty]}|\varphi'(y)|
	\\
	&\leq 2 \gamma^2 \sup_{y \in [-||X||_\infty, ||X||_\infty]}|\varphi'(y)|^2+ 2\sup_{\varepsilon < \varepsilon_0}  \sup_{t \in [0,\,T]}[X,X]_t^\varepsilon\, \delta(\varphi, \,3 \gamma) \sup_{y \in [-||X||_\infty, ||X||_\infty]}|\varphi'(y)|,
\end{align*}
where in the latter inequality we have used \eqref{estXsquare}. 
Since $\gamma$ is arbitrary and $\varphi'$ is uniformly continuous on compact intervals, then 
\begin{align*}
	J_1(\cdot, \varepsilon) \underset{\varepsilon \rightarrow 0}{\rightarrow} \sum_{s \leq \cdot}\Big[\Big(\int_0^1 \varphi'(X_{s-} + a \Delta X_s) da\Big)^2 - (\varphi'(X_{s-}))^2\Big](\Delta X_s)^2, 
\end{align*}
uniformly in $t$ for the fixed $\omega$. In particular, this implies \eqref{toproveJ1}.

By  \eqref{J2conv} and \eqref{toproveJ1}, and the fact that $[X, X]= [X, X]^c + \sum_{s \leq t}(\Delta X_s)^2$,  \eqref{J1+J2} yields, for $t \in [0,T]$,  
\begin{align*}
	&\frac{1}{\varepsilon}\int_0^t(\varphi(X_{(s + \varepsilon)\wedge t})-\varphi(X_{s}))^2 ds\\
	&\underset{\varepsilon \rightarrow 0}{\longrightarrow}
	\int_0^t (\varphi'(X_{s-}))^2 d [X,X]^c_s+\sum_{s \leq t}(\varphi'(X_{s-}))^2(\Delta X_s)^2\\
	&\qquad +\sum_{s \leq t}\Big[\Big(\int_0^1 \varphi'(X_{s-} + a \Delta X_s) da\Big)^2 - (\varphi'(X_{s-}))^2\Big](\Delta X_s)^2\\
	&=\int_0^t (\varphi'(X_{s-}))^2 d [X,X]^c_s+\sum_{s \leq t}\Big(\int_0^1 \varphi'(X_{s-} + a \Delta X_s) da\Big)^2(\Delta X_s)^2, \quad \textup{u.c.p.}
\end{align*}
The result follows because 
$$
\Delta \varphi(X_s) = \varphi(X_s)-\varphi(X_{s-})= \Delta X_s\int_0^1 \varphi'(X_{s-} + a \Delta X_s) da.
$$

\noindent  2. The result follows from point 1 by polarity arguments. 
\endproof



\section{Back to weak Dirichlet processes}
\label{S:WDP}

Let  $(\mathcal{F}_t)_{t \in [0,\,T]}$ be a filtration fulfilling the usual conditions on $(\Omega, \mathcal F, \P)$.
Let $X = (X_t)_{t \in [0,T]}$ be an $(\mathcal F_t)$-weak Dirichlet process with finite quadratic  variation.
 Cutting large jumps via a truncation function $k$, by Corollary 3.17 in \cite{BandiniRusso_RevisedWeakDir}, 
	the process
	\begin{equation}\label{specX}
	X_t - \sum_{s \le t} (\Delta X_s - k(\Delta X_s))=X_t - (x-k(x))\star \mu^X_t,\quad  t \in [0,T],
	\end{equation}
	is a special weak Dirichlet process, and therefore it admits a unique decomposition 
	$$
	X - (x-k(x))\star \mu^X =  B^{k,X}+M^{k,X}, 
	$$
	with $M^{k,X}$ an  $(\mathcal F_t)$-martingale and $B^{k,X}$ an $(\mathcal F_t)$-martingale orthogonal and predictable process. By Remark 3.22 in \cite{BandiniRusso_RevisedWeakDir}, this decomposition has 
	 the form
$$
X - (x-k(x))\star \mu^X= B^{k, X} + X^c + k(x) \star (\mu^X - \nu^X),
$$  
where $X^c$ is the unique local martingale component of $X$, $\nu^X$ is the compensator of $\mu^X$.
In particular,  
$
X = M+A
$
 with 
\begin{align}
	M&:= X^c + k(x) \star (\mu^X -\nu^X),\label{Mbis}\\
	A&:= B^{k,X} + (x-k(x))\star \mu^X. \label{Abis}
\end{align}
\begin{remark}\label{R:34}
The unique decomposition of  the  weak Dirichlet process $X$, provided in  Proposition 3.2 in \cite{BandiniRusso_RevisedWeakDir},  is therefore
\begin{equation}\label{uniqdecbis}
	X= X^c + \Gamma
\end{equation}
with 
$$
\Gamma = k(x) \star (\mu^X -\nu^X)+B^{k,X}  + (x-k(x))\star \mu^X.
$$
We also define $C^X:= [X^c, X^c]$. The triplet $(B^{k,X}, C^X, \nu^X)$ will be associated to the characteristics of $X$ when it is an $(\mathcal F_t^X)$-weak Dirichlet process, see Section \ref{chp4} below.
\end{remark}
	By Remark \ref{R:34} we easily see that 
	\begin{equation}\label{R3.1}
	\Delta B^{k,X}_t = \int_\R k(x) \nu^X(\{t\}\times dx). 
	\end{equation}
	This property is classical for semimartingales, see formula (2.14), Section II, in \cite{JacodBook}. We recall that when $X$ is a semimartingale, $B^{k,X}$ has bounded variation. This is no longer the case in the present framework. Nevertheless, we have the following. 


\begin{theorem}
\label{L:Bfinitevar}
	Let $X$ be an $(\mathcal F_t)$-weak Dirichlet process with finite quadratic  variation with decomposition \eqref{uniqdecbis}. 
	Then $B^{k,X}$ is a finite quadratic variation process, with  
  \begin{align}\label{brack_B2_bis}
 	[B^{k,X}, B^{k,X}] 
 	&=[X,X]^c -[X^c,X^c] 
 	+\sum_{s \leq \cdot} \Big|\int_\R k(x)  \nu^X(\{s\} \times dx)\Big|^2. 
 \end{align} 
\end{theorem}

\proof 
Recalling that  
$
X = M+A
$
 with $M$ and $A$ provided respectively by \eqref{Mbis} and \eqref{Abis}, we have 
 \begin{align*}
 	B^{k,X}= X- X^c - k(x) \star (\mu^X -\nu^X)- (x-k(x))\star \mu^X. 
 \end{align*} 
 From now on, let $t \in [0,T]$. 
 Since bounded variation processes and purely discontinuous martingales are martingale orthogonal processes, it follows that
 \begin{align}\label{brack_B}
 	[B^{k,X}, B^{k,X}]_t &= [X,X]_t +[X^c,X^c]_t +[k(x) \star (\mu^X -\nu^X),k(x) \star (\mu^X -\nu^X)]_t\notag\\
 	&+[(x-k(x))\star \mu^X,(x-k(x))\star \mu^X]_t-2[X,X^c]_t-2 [X,k(x) \star (\mu^X -\nu^X)]_t\notag\\
 	&- 2[X, (x-k(x))\star \mu^X]_t+2[k(x) \star (\mu^X -\nu^X), (x-k(x))\star \mu^X]_t,
 \end{align}
 provided  the right-hand side is well-defined.  
Now we notice that, by \eqref{uniqdecbis}, $[X, X^c]= [X^c, X^c]$. Moreover,  since $(x-k(x))\star \mu^X$ 
is a bounded variation process, by Proposition 2.14 in \cite{BandiniRusso1},
\begin{align*}
[(x-k(x))\star \mu^X,(x-k(x))\star \mu^X]_t&= \sum_{s \leq t} |\Delta X_s-k(\Delta X_s)|^2, \\
[X, (x-k(x))\star \mu^X]_t&=\sum_{s \leq t}\Delta X (\Delta X_s-k(\Delta X_s)),\\
 [k(x) \star (\mu^X -\nu^X), (x-k(x))\star \mu^X]_t&= \sum_{s \leq t}(\Delta X_s-k(\Delta X_s))\int_\R k(x)  (\mu^X -\nu^X)(\{s\} \times dx).
\end{align*}
On the other hand, being $k(x) \star (\mu^X -\nu^X)$ a purely discontinuous martingale, by Proposition 5.3 in \cite{BandiniRusso1},
$$
[k(x) \star (\mu^X -\nu^X),k(x) \star (\mu^X -\nu^X)]_t= \sum_{s \leq t} \Big|\int_\R k(x)  (\mu^X -\nu^X)(\{s\} \times dx)\Big|^2.
$$
 Finally, by Lemma \ref{L:YZ} with  $Y= k(x) \star (\mu^X -\nu^X)$ and $Z= X$,  
\begin{align*}
	[X,k(x) \star (\mu^X -\nu^X)]_t&=\sum_{s \leq t}\Delta X_s \int_\R k(x)  (\mu^X -\nu^X)(\{s\} \times dx).
\end{align*}
Plugging previous terms in  \eqref{brack_B}  we get  
\begin{align}\label{brack_B2}
 	[B^{k,X}, B^{k,X}]_t 
 	&=[X,X]_t -[X^c,X^c]_t + \sum_{s \leq t} \Big|\int_\R k(x)  (\mu^X -\nu^X)(\{s\} \times dx)\Big|^2\notag\\
 	&+\sum_{s \leq t} |\Delta X_s-k(\Delta X_s)|^2-2 \sum_{s \leq t}\Delta X_s \int_\R k(x)  (\mu^X -\nu^X)(\{s\} \times dx)\notag\\
 	&-2\sum_{s \leq t}\Delta X (\Delta X_s-k(\Delta X_s))+2\sum_{s \leq t}(\Delta X_s-k(\Delta X_s))\int_\R k(x)  (\mu^X -\nu^X)(\{s\} \times dx)\notag\\
 	&=[X,X]_t -[X^c,X^c]_t + \sum_{s \leq t} \Big|\int_\R k(x)  (\mu^X -\nu^X)(\{s\} \times dx)\Big|^2\notag\\
 	&+\sum_{s \leq t} |\Delta X_s-k(\Delta X_s)|^2-2\sum_{s \leq t}\Delta X (\Delta X_s-k(\Delta X_s))\notag\\
 	&-2\sum_{s \leq t}k(\Delta X_s)\int_\R k(x)  (\mu^X -\nu^X)(\{s\} \times dx), 
 \end{align} 
 which implies in particular that $[B^{k,X}, B^{k,X}]_t$ is finite. 
 Now, recalling that $[X,X]_t = [X,X]^c_t + \sum_{s \leq t}|\Delta X_s|^2$, and noticing that 
\begin{align*}
 \sum_{s \leq t} |\Delta X_s-k(\Delta X_s)|^2&= \sum_{s \leq t} |\Delta X_s|^2+\sum_{s \leq t} |k(\Delta X_s)|^2-2 \sum_{s \leq t}\Delta X_s \, k(\Delta X_s),\\
-2 \sum_{s \leq t}\Delta X (\Delta X_s-k(\Delta X_s))&= -2\sum_{s \leq t} |\Delta X_s|^2 +2\sum_{s \leq t}\Delta X_s \, k(\Delta X_s), 
\end{align*}
  formula  \eqref{brack_B2} reads 
 \begin{align}\label{brack_B2_tris}
 	[B^{k,X}, B^{k,X}]_t 
 	&=[X,X]_t -[X^c,X^c]_t  -\sum_{s \leq t} |\Delta X_s|^2 +\sum_{s \leq t} |k(\Delta X_s)|^2
 	\\
 	&+ \sum_{s \leq t} \left(\Big|\int_\R k(x)  (\mu^X -\nu^X)(\{s\} \times dx)\Big|^2 -2 \,k(\Delta X_s) \int_\R k(x)  (\mu^X -\nu^X)(\{s\} \times dx)\right). \notag
 \end{align} 
 Finally, we notice that 
  \begin{align}\label{toplugfinal}
 	&\Big|\int_\R k(x)  (\mu^X -\nu^X)(\{s\} \times dx)\Big|^2-2 k(\Delta X_s) \int_\R k(x)  (\mu^X -\nu^X)(\{s\} \times dx)\notag\\
 	&=|k(\Delta X_s)|^2+ \Big|\int_\R k(x)  \nu^X(\{s\} \times dx)\Big|^2-2 \,k(\Delta X_s) \int_\R k(x)  \nu^X(\{s\} \times dx)\notag\\
 	&-2 |k(\Delta X_s)|^2+2 k(\Delta X_s)\int_\R k(x)  \nu^X(\{s\} \times dx)\notag\\
 	&= \Big|\int_\R k(x)  \nu^X(\{s\} \times dx)\Big|^2- |k(\Delta X_s)|^2. 
 \end{align}
 Plugging \eqref{toplugfinal} in \eqref{brack_B2_tris} (noticing that $[X,X]_t -\sum_{s \leq t} |\Delta X_s|^2= [X,X]^c_t$) we get 
  \eqref{brack_B2_bis}.
\endproof

 \begin{corollary}\label{C_bracketB}
 	Let $X$ be an $(\mathcal F_t)$-weak Dirichlet process with finite quadratic  variation with decomposition \eqref{uniqdecbis}. Then 
 		  \begin{align}\label{33}
 	[B^{k,X}, B^{k,X}] 
 	&=[X,X]^c -[X^c,X^c] 
 	+\sum_{s \leq \cdot} |\Delta B^{k,X}_s|^2, 
 \end{align}
  	or equivalently,
 	\begin{equation}\label{33bis}
 	[X, X]^c =  [X^c, X^c] + [B^{k,X},B^{k,X}]^c. 
 	\end{equation}
 \end{corollary}
\proof  
Taking into  account \eqref{R3.1}, formula \eqref{brack_B2_bis} of Theorem \ref{L:Bfinitevar} can be rewritten as \eqref{33}, which is in turn 
 equivalent to \eqref{33bis}.
\endproof
 \begin{remark}\label{R:3.5}
If $X$ is a semimartingale (resp. a Dirichlet process), by Corollary \ref{C_bracketB} we recover the result 
	\begin{equation}\label{result}
	[X,X]^c =[X^c,X^c],  
	\end{equation}
	proved in Proposition 3.4 in \cite{BandiniRusso_RevisedWeakDir} (resp. in Proposition 6.2 in \cite{BandiniRusso_DistrDrift}). 
	In fact we have the following. 
	\begin{itemize}
\item[(i)]
Let $X$ be a semimartingale. Then $B^{k,X}$  is a finite variation process, so by Proposition 3.14 in \cite{BandiniRusso1}
	$$
	[B^{k,X}, B^{k,X}]_t= \sum_{ s\leq t} |\Delta B_s^{k,X}|^2,\quad t \in  [0,T], 
	$$
	 and therefore by \eqref{33} we recover \eqref{result}. 
	 \item[(ii)] Let $X$ be a Dirichlet process. Then 
	 it is a special weak Dirichlet process, and by \eqref{Mbis}-\eqref{Abis}, taking into account Corollary 3.21-(ii) in \cite{BandiniRusso_RevisedWeakDir}, it admits the unique decomposition 
	 $$
	 X= X^c + x \star (\mu^X-\nu^X) + (x-k(x)) \star \nu^X  + B^{k,X}. 
	 $$
	 Let $Y= (x-k(x)) \star \nu^X$.
	 Since $X$ is a Dirichlet process, 
	 $$
	 [ Y  + B^{k,X},  Y  + B^{k,X}] = 0, 
	 $$
	 so that 
	 $$
	 [B^{k,X}, B^{k,X}] = -[Y,Y] - 2 [Y, B^{k,X}]= - \sum_{s \leq \cdot} |\Delta Y_s|^2 - 2 [Y, B^{k,X}].
	 $$
	 On the other hand, by Lemma 2.4 	 
$$
[Y,B^{k,X}]= \sum_{s \leq \cdot} \Delta Y_s \Delta B_s^{k,X}, 
$$
so that 
$$
 [B^{k,X}, B^{k,X}]_t = - \sum_{s \leq t} |\Delta Y_s|^2 - 2 \sum_{s \leq t} \Delta Y_s \Delta B_s^{k,X}.
	 $$
	 We get that $[B^{k,X}, B^{k,X}]^c_t =0$,
	 so by \eqref{33bis} we recover \eqref{result}. 
\end{itemize}
	 
\end{remark}


We state here a slight modification of Theorem 3.37 in \cite{BandiniRusso_RevisedWeakDir}.
\begin{theorem}\label{T:new3.36}
	Let $X$ be an $(\mathcal F_t)$-weak Dirichlet process with finite quadratic variation, taking values in an open interval $\mathcal O$. Let $v\in C^{0,1}([0,T] \times \mathcal O)$.  Then  $Y_t = v(t, X_t)$ is an $(\mathcal F_t)$--weak Dirichlet with continuous martingale component 
	\begin{align}\label{Yc}
	Y^c = Y_0 + \int_{0}^{\cdot}\partial_x v(s, X_{s})\,dX^c_s.
\end{align}
\end{theorem}
\proof
The proof follows the same lines of the one of Theorem 3.37 in \cite{BandiniRusso_RevisedWeakDir}, taking into account that the set $\mathcal O$ is open and convex.

\section{Main result}\label{chp4}
\subsection{Characteristics of a weak Dirichlet process}
We denote by $\check \Omega$ the canonical  space of all c\`adl\`ag functions $\check \omega: [0,T] \rightarrow \R$,  and   by $\check X$ the canonical process defined by  $\check X_t(\omega)= \check \omega(t)$. We also set $\check {\mathcal  F}= \sigma(\check X)$, and $\check {\mathbb  F}= (\check {\mathcal  F}_t)_{t \in [0,T]}$. 
We suppose below that $X$ is an $(\mathcal F_t^X)$-weak Dirichlet process.
Let $\mu$ be the jump measure of $\check X$ and $\nu$ the compensator of $\mu$ under the law $\mathcal L(X)$ of $X$. 
\begin{definition}\label{D:genchar}
We call \emph{characteristics} of $X$, associated with $k \in \mathcal K$, 
  the triplet $(B^k,C, \check X^c\rangle,\nu)=(B^{k,\check X},C^{\check X},\nu^{\check X})$ on $(\check \Omega, \check {\mathcal F}, \check {\mathbb F})$ obtained from the unique decomposition \eqref{uniqdecbis} in Remark \ref{R:34} for $\check X$ under $\mathcal L(X)$. 
  In particular, 
		 \begin{itemize}
	\item[(i)] $B^k$ is a predictable and $\check{\mathbb F}$-martingale orthogonal process, with $B_0^k=0$;
	\item[(ii)] $C$ is an  $\check{\mathbb F}$-predictable and  increasing process, with $C_0=0$; 	
	\item[(iii)] $\nu$ is an $\check{\mathbb F}$-predictable random measure on $[0,T]\times \R$.
	\end{itemize}
	\end{definition}

 Let $X = (X_t)_{t \in [0,T]}$ be an $(\mathcal F_t^X)$-weak Dirichlet process with finite quadratic  variation with characteristics $(B^k, C, \nu)$. By Remark 3.26 in \cite{BandiniRusso_RevisedWeakDir},  
 we have $B^{k,X}= B^k \circ X$, $\nu^X= \nu \circ X$ and  $C^{X}= C \circ X$.
 Therefore in this case \eqref{Mbis}-\eqref{Abis}
read
\begin{align}
	M&:= X^c + k(x) \star (\mu^X -(\nu \circ X)),\label{M}\\
	A&:= B^{k}\circ X + (x-k(x))\star \mu^X.\label{A}
\end{align}

\subsection{The equivalence theorem}

We  provide an equivalence result for c\`adl\`ag weak Dirichlet processes that extends the analogous one for c\`adl\`ag semimartingales, see Theorem \ref{T:equiv_mtgpb_semimart}. 

\begin{theorem}\label{T:to_prove}
Let  $X$ be a  c\`adl\`ag  process with  finite quadratic variation. 
Let $B^k$ be an  $({\mathcal {\check F}_t})$-predictable process  such that $B_0^k=0$ and  $B^k \circ X$ has finite quadratic variation, 
 $C$ be
	an $({\mathcal {\check F}_t})$-adapted continuous process  such that  $C_0=0$ and  $C \circ X$ has finite variation, 	and  $\nu$ be an $({\mathcal {\check F}_t})$-predictable random measure on $[0,T] \times \R$.
		
There is equivalence between the two following statements.
\begin{itemize} \item [(i)] X is an $({\mathcal F_t^X})$-weak Dirichlet process with local characteristics $(B^k, C, \nu)$.
\item [(ii)] For each bounded function $F$ of class $C^{1,2}$, 
$
\int_{]0,\cdot]}  \partial_x F(s,X_{s}) \,d^- (B^{k}\circ X)_s
$
is an $({\mathcal F_t^X})$-martingale orthogonal process,
 and the process 
\begin{align}\label{MF}
&F(\cdot, X_{\cdot}) - F(0, X_0) -\int_0^\cdot \partial_s F(s, X_s) ds- \frac{1}{2} \int_0^{\cdot} \partial^2_{xx} F(s,X_s) \,(d(C \circ X)_s +d[B^k \circ X, B^k \circ X]^c_s)\notag \\
&-  \int_{]0,\cdot]}  \partial_x F(s,X_{s}) \,d^{-} (B^{k}\circ X)_s
   -(F(\cdot,X_{-} + x) -F(\cdot,X_{-})-k(x)\,\partial_xF(\cdot,X_{-}))\star(\nu\circ X)
\end{align} 
is an $({\mathcal F_t^X})$-local martingale. 
\end{itemize} 
\end{theorem}

\begin{remark}
\begin{enumerate}
\item
  The stochastic integral $\int_{]0,\cdot]}  \partial_x F(s,X_{s}) \,d^{-} (B^{k}\circ X)_s$	is a predictable process.
Indeed,  its jump process is given by 
\begin{equation} \label{JumpPred}
\partial_x F(t,X_{t-}) \,\Delta  (B^{k}\circ X)_t, \quad t \in [0,T],
\end{equation}
see  $(1.15)_{-}$  in \cite{rv95}.
Now, the first term of the product in \eqref{JumpPred}  is a c\`agl\`ad,
therefore predictable, the second one is also predictable since $B^k \circ X$ is predictable. 
\item  If $X$ is a semimartingale, then $B^{k,X}$ has bounded variation, so that 
$$\int_{]0,\cdot]}  \partial_x F(s,X_{s}) \,d^{-} (B^{k}\circ X)_s=\int_0^{\cdot}  \partial_x F(s,X_{s-}) \,d (B^{k}\circ X)_s$$ has bounded variation, therefore it is an $(\mathcal F_t)$-martingale orthogonal process. Moreover, $[B^k \circ X, B^k \circ X]^c_s=0$, see Remark \ref{R:3.5}-(i), and we retrieve the result of Jacod-Shiryaev, see Theorem \ref{T:equiv_mtgpb_semimart} (and \eqref{JacodEquiv} in the non-homogeneous form). 

\end{enumerate}
\end{remark}

\proof $(i) \Rightarrow (ii)$.
Let $X$ be an $(\mathcal F_t^X)$-weak Dirichlet process with finite quadratic  variation with characteristics $(B^k, C, \nu)$. Let $F\in C^{1,2}$.
We recall that $X = M + A$ as in \eqref{M} and \eqref{A}.
By  Theorem 5.15  in \cite{BandiniRusso1} 
we have \begin{align}\label{CO1_Ito_formula_weak_D} 
 F(t,X_t)&=F(0,X_0)+\int_0^t\partial_x F(s,X_{s-})\,d M_s\nonumber\\ 
&+  (F(\cdot,X_{-}+x)-F(\cdot,X_{-}))\frac{k(x)}{x}\star (\mu^X-\nu^X)_t\nonumber\\
&-  x\,\partial_x F(\cdot,X_{-})\,\frac{k(x)}{x}\star (\mu^X-\nu^X)_t\nonumber\\ 
&+ (F(\cdot,X_{-}+x)-F(\cdot,X_{-})-x\,\partial_x F(\cdot,X_{-}))\,\frac{x-k(x)}{x}\star \mu^X_t\notag\\
&+\int_0^t\partial_s F(s,X_s)\,ds + \int_{]0,t]}\partial_x F(s,X_{s})\,d^{-}A_s +\frac{1}{2} \int_0^t\partial_{xx}^2 F(s,X_s)\,d[X,X]_s^c\nonumber\\ 
&  + 
 (F(\cdot,X_{-}+x)-F(\cdot,X_{-})-x\,\partial_x F(\cdot ,X_{-}))\frac{k(x)}{x}\star \nu^X_t.
\end{align}
In fact Theorem 5.15  in \cite{BandiniRusso1} was written for $k(x) = x \one_{\{|x|\leq 1\}}$, but it naturally extends to a generic truncation function $k$.

By \eqref{A}  and using Proposition \ref{P:3.1},   we get 
\begin{align}\label{dA}
\int_{]0,t]}\partial_x F(s,X_{s})\,d^{-}A_s= \int_{]0,t]}\partial_x F(s,X_{s})\,d^{-} B^{k, X}_s 
+ \partial_x F(\cdot,X_{-})\,(x-k(x))\star \mu^X_t.
\end{align}
On the other hand, by \eqref{M} 
\begin{align*}
\int_0^t\partial_x F(s,X_{s-})\,dM_s= \int_0^t\partial_x F(s,X_{s-})\,d X^c_s + \int_0^t\partial_x F(s,X_{s-})\,d M^{d,k}_s, 
\end{align*}
where $M^{d,k}= k(x) \star (\mu^X -\nu^X)$. We notice that  
\begin{align*}
\Delta \left(\int_0^t\partial_x F(s,X_{s-})\,d M^{d,k}_s\right)&= \partial_x F(t,X_{t-})\,\Delta M^{d,k}_t\\
&= \partial_xF(t,X_{t-})\int_\R  k(x)(\mu^X-\nu^X)(\{t\} \times dx)\\
&= \Delta\left(\partial_x F(\cdot,X_{-})\,k(x)\star (\mu^X-\nu^X)_t\right).
\end{align*}
We remind that, for any $Y(\cdot)$ predictable random field,  $ Y(x)  \star (\mu^X-\nu^X)$ is the unique purely discontinuous martingale orthogonal process whose jumps are indistinguishable from 
$$ 
\Delta\left( Y(x)\star (\mu^X-\nu^X)\right), 
$$
see  Corollary 4.19, Section I, in \cite{JacodBook}.
We conclude that 
\begin{align}\label{dM}
&\int_0^t\partial_x F(s,X_{s-})\,dM_s
= \int_0^t\partial_x F(s,X_{s-})\,d X^c_s + 
\partial_x F(\cdot,X_{-})\,k(x)\star (\mu^X-\nu^X)_t.  
\end{align}
Plugging  \eqref{dA} and \eqref{dM} in \eqref{CO1_Ito_formula_weak_D}, we get
\begin{align} F(t,X_t)&=F(0,X_0)
+\int_0^t\partial_x F(s,X_{s-})\,d X^c_s +\partial_x F(\cdot,X_{-})\,k(x)\star  (\mu^X -\nu^X)_t\nonumber\\ 
&+ (F(\cdot,X_{-}+x)-F(\cdot,X_{-}))\,\frac{k(x)}{x}\star(\mu^X-\nu^X)_t\nonumber\\ 
&- \partial_x F(\cdot,X_{-})\,k(x)\star (\mu^X-\nu^X)_t\nonumber\\ 
&+ (F(\cdot,X_{-}+x)-F(\cdot,X_{-})-x\,\partial_x F(\cdot,X_{-}))\,\frac{x-k(x)}{x}\star\mu^X_t\notag\\
&+\int_0^t\partial_s F(s,X_s)\,ds  +\frac{1}{2} \int_0^t\partial_{xx}^2 F(s,X_s)\,d[X,X]_s^c\nonumber\\ 
&  +  (F(\cdot,X_{-}+x)-F(\cdot,X_{-})-x\,\partial_x F(\cdot,X_{-}))\frac{k(x)}{x}\star\nu^X_t\notag\\
&+ \int_{]0,t]} \partial_x F(s,X_{s})\,d^{-}B^{k,X}_s
+ \partial_x F(\cdot,X_{-})\,(x-k(x)) \star\mu^X_t,\notag
\end{align}
that reads
\begin{align}\label{CO1_Ito_formula_weak_D_BIS} 
 F(t,X_t)&=F(0,X_0)+\int_0^t\partial_x F(s,X_{s-})\,d X^c_s\nonumber\\ 
&+ (F(\cdot,X_{-}+x)-F(\cdot,X_{-}))\,\frac{k(x)}{x}\star(\mu^X-\nu^X)_t\nonumber\\ 
&+\int_0^t\partial_s F(s,X_s)\,ds  +\frac{1}{2} \int_0^t\partial_{xx}^2 F(s,X_s)\,d[X,X]_s^c+ \int_{]0,t]}\partial_x F(s,X_{s})\,d^{-}B^{k,X}_s\nonumber\\ 
&  + (F(\cdot,X_{-}+x)-F(\cdot,X_{-})-x\,\partial_x F(\cdot,X_{-}))\frac{k(x)}{x}\star \nu^X_t\notag\\
&+(F(\cdot,X_{-}+x)-F(\cdot,X_{-}))\,\frac{x-k(x)}{x}\star\mu^X_t.
\end{align}
At this point we make use of the fact that $F$ is bounded. Indeed, being $F \in C_b^{0,1}$, by Theorem  3.15 and Remark 3.16 in \cite{BandiniRusso_RevisedWeakDir},  	\begin{equation}\label{intY}
\textup{$\forall \,a \in \R_+$ s.t.}\,\,\,	(F(\cdot,X_{-}+x)-F(\cdot,X_{-})) \,\one_{\{|x| >a\}} \star \mu^X \in \mathcal{A}_{\textup{loc}}.
	\end{equation}
Therefore, by Lemma C.3 in \cite{BandiniRusso_RevisedWeakDir},  
$$
(F(\cdot,X_{-}+x)-F(\cdot,X_{-})) \,\frac{x-k(x)}{x} \star \mu^X \in \mathcal{A}_{\textup{loc}},
$$
so that 
$$
(F(\cdot,X_{-}+x)-F(\cdot,X_{-})) \,\frac{x-k(x)}{x} \star \nu^X_t
$$
 is well-defined for every $t \in \R_+$. Adding and subtracting the above mentioned term in \eqref{CO1_Ito_formula_weak_D_BIS}, and using Corollary \ref{C_bracketB} (recalling that $C^X= [X^c,X^c]$),  we get 
\begin{align}\label{CO1_Ito_formula_weak_D_TRIS} 
 F(t,X_t)&=F(0,X_0)+\int_0^t\partial_x F(s,X_{s-})\,d X^c_s
+ (F(\cdot,X_{-}+x)-F(\cdot,X_{-}))\star(\mu^X-\nu^X)_t\nonumber\\   
&+\int_0^t\partial_s F(s,X_s)\,ds  +\frac{1}{2} \int_0^t\partial_{xx}^2 F(s,X_s)\,(d C^X_s + d [B^{k,X}, B^{k,X}]^c_s) \notag\\
&+ \int_{]0,t]}\partial_x F(s,X_{s})\,d^{-}B^{k,X}_s\notag\\
& +  (F(\cdot,X_{-}+x)-F(\cdot,X_{-})-k(x)\,\partial_x F(\cdot,X_{-}))\star\nu^X_t,
\end{align}
which implies that \eqref{MF} is an $(\mathcal F^X_t)$-local martingale.

It remains to prove that  $\int_{]0,\cdot]}\partial_x F(s,X_{s})\,d^{-}B^{k,X}_s$
is martingale orthogonal. 
By Theorem 3.37 in \cite{BandiniRusso_RevisedWeakDir},  $F(t, X_t)$ is an  $(\mathcal F_t^X)$-weak Dirichlet process
 with  continuous martingale component  
$$
F(0,X_0) + \int_0^\cdot \partial_x F(s,X_{s-})\,d X^c_s.
$$
By the uniqueness of the decomposition of weak Dirichlet processes (see Proposition 3.2 in \cite{BandiniRusso_RevisedWeakDir}), we get from \eqref{CO1_Ito_formula_weak_D_TRIS}  that 
\begin{align} \label{OrthPart}
 \int_{]0,\cdot]} \partial_x F(s,X_{s})\,d^{-}B^{k,X}_s + \Gamma
 \end{align}
 is a martingale orthogonal process, where
\begin{align*}
\Gamma:=& 
(F(\cdot,X_{-}+x)-F(\cdot,X_{-}))\star (\mu^X-\nu^X)\nonumber\\  
&+\int_0^\cdot \partial_s F(s,X_s)\,ds  +\frac{1}{2} \int_0^t\partial_{xx}^2 F(s,X_s)\,(d C^X_s + d [B^{k,X}, B^{k,X}]^c_s) \notag\\
& +  (F(\cdot,X_{-}+x)-F(\cdot,X_{-})-k(x)\,\partial_x F(\cdot,X_{-}))\star \nu^X. 
\end{align*}
Here $\Gamma$ is a martingale orthogonal process since is the sum of a purely discontinuous  martingale and bounded variation processes.
 This finally allows to conclude that
 the first term in \eqref{OrthPart} is martingale orthogonal.

\medskip 
\noindent 
$(ii) \Rightarrow (i)$. We apply (ii) with $F$ time-homogeneous. Then, for each  function $F$ of class $C^{2}_b$, 
$$
\int_{]0,\cdot]}  F'(X_{s}) \,d^- (B^{k}\circ X)_s
$$
is an $({\mathcal F_t})$-martingale orthogonal process, 
 and the process 
\begin{align}\label{C1}
&M^F := F( X_{\cdot}) - F( X_0) - \frac{1}{2} \int_0^{\cdot}  F''(X_s) \,(d(C \circ X)_s +d[B^{k}\circ X, B^{k}\circ X]^c_s)\notag\\
&
-  \int_{]0,\cdot]}   F'(X_{s}) \,d^{-} (B^{k}\circ X)_s
-  (F(X_{-} + x) -F(X_{-})-k(x)\,F'(X_{-}))\star (\nu\circ X), 
\end{align} 
is an $(\mathcal F_t^X)$-local martingale. 


\medskip 

\noindent \emph{Step 1: $X$ is an $({\mathcal F^X_t})$-weak Dirichlet process.} 
Let us now set $\tilde F(x)= \arctan{x}$. Being $\tilde F \in C^{2}_b$, by Theorem 3.15 and Remark 3.16 in \cite{BandiniRusso_RevisedWeakDir},   $Y = \tilde F(X)$ is an $(\mathcal F_t^X)$-special weak Dirichlet process. On the other hand, by  Lemma \ref{L:app1},  $Y$ has finite quadratic variation. Since $X_s= \tilde F^{-1}(Y_s)$ and $\tilde F^{-1} \in C^{1}(\mathcal O)$ with
$\mathcal O=]-\frac{\pi}{2}, \frac{\pi}{2}[$, we can apply Theorem \ref{T:new3.36} to $X_s = \tilde F^{-1}(Y_s)$,  getting that $X$ is  a weak Dirichlet process.

\medskip 

\noindent \emph{Step 2: if $(B^k, C,   \nu)$ and $(\bar B^k, \bar C, \bar\nu)$ both verify  \eqref{C1}, then 
\begin{equation}\label{AimStep2}
	(B^k, C,   \nu)=(\bar B^k, \bar C, \bar\nu)
\end{equation}
in the $\mathcal L(X)$-sense. }


 
 \noindent \emph{Step 2a:  
 $B^{k}-\bar B^{k}$ has finite variation.
  }
We start by noticing that, for every $F \in C^2_b$, by  uniqueness of decomposition of the special weak Dirichlet process $F(X)$, by \eqref{C1} we have 
  \begin{align*}
 &\frac{1}{2} \int_0^{\cdot} F''(X_s) \,(d(C\circ X-\bar C\circ X)_s +d([B^{k}\circ X, B^{k}\circ X]^c-[\bar B^{k}\circ X, \bar B^{k}\circ X]^c)_s)\notag\\
 &+ \int_{]0,\cdot]}  F'( X_{s}) \,d^{-} (B^{k}\circ X-\bar B^{k}\circ X)_s\notag\\
&+   (F( X_{-} + x) -F( X_{-})-k(x)\,F'( X_{-}))\star (\nu\circ X-\bar \nu\circ X)=0. 
\end{align*}
By Remark 3.24-ii) in \cite{BandiniRusso_RevisedWeakDir}, $[ B^{k}\circ X,  B^{k}\circ X]= [B^k, B^k]\circ X$ almost surely with respect to $\mathcal L(X)$. So, again almost surely with respect to $\mathcal L(X)$,  $[ B^{k}\circ X,  B^{k}\circ X]^c= [B^k, B^k]^c\circ X$, since $\Delta ( B^{k}\circ X)= (\Delta B^k)\circ X$, taking into account \eqref{QVC}.  Therefore, 
 \begin{align}\label{2}
 &\frac{1}{2} \int_0^{\cdot} F''(\check X_s) \,(d(C-\bar C)_s +d([B^{k}, B^{k}]^c-[\bar B^{k}, \bar B^{k}]^c)_s)+ \int_{]0,\cdot]}  F'(\check X_{s}) \,d^{-} (B^{k}-\bar B^{k})_s\notag\\
&+   (F(\check X_{-} + x) -F(\check X_{-})-k(x)\,F'(\check X_{-}))\star (\nu-\bar \nu)=0. 
\end{align}
It will enough to prove that $B^{k}-\bar B^{k}$ has finite variation on 
$$
\check \Omega_n := \{\omega \in \check \Omega: \,\,\sup_{t \leq T} |\check X_t| \leq n\},
$$
 being $\check \Omega = \cup_n\check \Omega_n$ up to a null set, with respect to $\mathcal L(X)$. We apply \eqref{2} with 
 $$
 F_n(x) =x \,\chi_n(x),
 $$
  where $\chi_n \in C^\infty_b$, $|\chi_n|\leq 1$ and 
$$
\chi_n(x)=
\begin{cases}
	1 & \textup{if} \quad |x|\leq  n,\\
	0 & \textup{if} \quad |x|>  n+1,
\end{cases}
$$
getting 
 \begin{align}\label{2bis}
 &\frac{1}{2} \int_0^{\cdot} F''_n(\check X_s) \,(d(C-\bar C)_s +d([B^{k}, B^{k}]^c-[\bar B^{k}, \bar B^{k}]^c)_s)+ \int_{]0,\cdot]} F'_n(\check X_{s}) \,d^{-} (B^{k}-\bar B^{k})_s\notag\\
&+  (F_n(\check X_{-} + x) -F_n(\check X_{-})-k(x)\,F'_n(\check X_{-}))\star (\nu-\bar \nu)=0. 
\end{align}
 By Remark \ref{R:12}, 
 \begin{align*}
 	&\one_{\check \Omega_n} \int_{]0,\cdot]}F'_n(\check X_s) d^- (B^k-\bar B^k)_s =  \one_{\Omega_n} \int_{]0,\cdot]}  d^- (B^k-\bar B^k)_s = \one_{\Omega_N} (B^k_t-\bar B^k_t);
 \end{align*}
 moreover, 
  	$$
  	\one_{\check \Omega_n} \int_0^\cdot F''_n(\check X_s) (d(C-\bar C)_s +d([B^{k}, B^{k}]^c-[\bar B^{k}, \bar B^{k}]^c)_s)=  0. 
  	$$
 Consequently, \eqref{2bis} on $\check \Omega_n$ yields  
  \begin{align*}
\one_{\check \Omega_n} (B^k-\bar B^k)=
- \one_{\check \Omega_n} (F_n(\check X_{-} + x) -F_n(\check X_{-})-k(x)\,F'_n(\check X_{-}))\star (\nu-\bar \nu),
\end{align*}
which implies that $\one_{\Omega_n} (B^k-\bar B^k)$ has finite variation. 

  \medskip 
 
 \noindent \emph{Step 2b: \eqref{AimStep2} holds true. }
 Let $u \in \R$.  Notice that  $(e^{iux} -1 - iuk(x)) \star  \nu \in \mathcal A_{\textup{loc}}$. Indeed,   $|e^{iux} -1 - iuk(x)|\leq \alpha (1 \wedge |x|^2)$ for some constant $\alpha$, and  $\sum_{s \leq \cdot}|\Delta X_s|^2 < +\infty $ a.s.,  hence $(1 \wedge |x|^2) \star \nu  \in \mathcal A_{\textup{loc}}$, see Proposition C.1 in \cite{BandiniRusso_RevisedWeakDir}. 

We remark that \eqref{2} extends for a complex valued function $F$ such that $Re(F)$ and $Im(F)$  belong to $C^2_b$. 
We can then apply \eqref{2} with $F(x) = e^{iux}$. We have  $F'(x)=iu F(x)$, $F''(x)=-u^2 F(x)$, and 
 $$
 F(\check X_{s-} +x) - F(\check X_{s-})-k(x)  F'(\check X_{s-}) = e^{iu \check X_{s-}}(e^{iux}-1 -iuk(x)).
 $$
 Since  by Step 2a the process $B^k - \bar B^k$ has bounded variation, by Remark \ref{R:1.1}-(i)
 $$
 \int_{]0,\cdot]} F'(\check X_s) d^- (B^k-\bar B^k)_s = \int_0^\cdot F'(\check X_{s-}) d(B^k-\bar B^k)_s, 
 $$
and 
 we get that
 \begin{align*}
& \int_{]0,\cdot]}  e^{iu \check X_{s-}}\Big[   iu  \,d(B^{k} -\bar B^{k})_s-\frac{1}{2}u^2 (d(C-\bar C)_s +d([ B^{k}, B^{k}]^c-[\bar B^{k}, \bar B^{k}]^c)_s)\\
&+\int_\R (e^{iux}-1 -iuk(x))(\nu-\bar \nu)(ds\,dx)\Big]=0 \quad \textup{up to an evanescent set},
\end{align*}
or,  equivalently,  for every $t \in [0,T]$, 
\begin{align}\label{ev}
\int_0^t  e^{iu \check X_{s-}}dH(u)_s=0\quad \textup{up to an evanescent set (with respect to $\mathcal L(X)$)},
\end{align}
with 
\begin{align}\label{ind}
 H(u)_t &:= iu (B^{k}-\bar B^k)_t - \frac{1}{2}u^2 ((C-\bar C)_t + ([B^{k}, B^{k}]^c-[\bar B^{k}, \bar B^{k}]^c)_t) \notag\\
 &+ \int_\R (e^{iux} -1 - iu \,k(x) )(\nu-\bar \nu)([0,t] \times dx ).
\end{align}  
Notice that, since  $B^{k}-\bar B^{k}$ has finite variation, the same property holds for $H(u)$. 

In particular, there is a null set $\mathcal N(u)$ for which \eqref{ev} holds for every $t$. 
Differentiating \eqref{ev} in $t$  for every $\omega \notin \mathcal N(u)$, we get that $H(u)_t=0$  for every $\omega \notin \mathcal N(u)$, for every $t$. 
We define $\mathcal N = \cup_{u \in \Q} \,\mathcal N(u)$.
Since the left-hand side of \eqref{ev} is continuous in $u$ (uniformly in $t$), $H(u)_t=0$  for every $\omega \notin \mathcal N$, for every $t$ and $u$.

Now let us fix $\omega \notin \mathcal N$ and $t \in [0,T]$. We set $b_t:= (B^k-\bar B^k)_t$, $c_t:= (C-\bar C)_t$,  $\Lambda_t(dx):= (\nu-\bar \nu)([0,t] \times dx)$, and $\bar b= \bar c= \bar \Lambda=0$. 
 Then, applying  Lemma \ref{L:app} below, 
we conclude that $b=0$, $c=0$ and $\Lambda=0$ and this concludes the proof of Step 2.   

 \medskip 
 
 \noindent \emph{Step 3: X is an $({\mathcal F_t^X})$-weak Dirichlet process with local characteristics $(B^k, C, \nu)$.} By Step 1, $X$ is an $({\mathcal F_t^X})$-weak Dirichlet process. Let $(\bar B^k, \bar C, \bar \nu)$ be the characteristics of $X$, see Definition \ref{D:genchar}. We conclude by applying $(i)\Rightarrow (ii)$ together with Step 2.
 \endproof
 
 The following lemma is the extension of Lemma 2.44, Chapter II, in \cite{JacodBook}   in the case of signed measures  $\Lambda$ and symmetric (not necessarily positive) matrices $c$. 
\begin{lemma}\label{L:app}
	Let $b, \bar b \in \R^d$, let $c, \bar c$  symmetric  $d \times d$ matrices, and $\Lambda, \bar\Lambda $  signed measures on $\R^d$ that satisfy $\Lambda(\{0\})=0$, $\bar \Lambda(\{0\})=0$ and whose total variation measure  integrate $(1 \wedge |x|^2)$.
	Let
	\begin{equation}\label{245}
		\psi(u) = iu  b - \frac{1}{2} u^T c  u + \int_{\R^d} (e^{iux}-1-iu k(x))\Lambda(dx), \quad u \in \R^d.
	\end{equation}
	If $\psi$ satisfies \eqref{245} with $(\bar b, \bar c, \bar \Lambda)$ also, then $b= \bar b$, $c = \bar c$ and $\Lambda= \bar \Lambda$. 
\end{lemma}

\proof 
Let $w \in \R^d \setminus \{0\}$ and define the function 
$$
\varphi_w(u):=\psi(u) -\frac{1}{2}\int_{-1}^1 \psi(u+s w)ds.
$$
One can easily prove that 
$$
\varphi_w(u)=\frac{1}{6} w^T c w + \int_{\R^d} \left(1- \frac{sin(wx)}{wx}\right)e^{iux} \Lambda(dx).
$$
Therefore, the function $\varphi_w(u)$ is the Fourier transform of the measure 
\begin{equation}\label{Gw}
G_w(dx)=\frac{1}{6} w^T c w \, \delta_0(dx) +  \left(1- \frac{sin(wx)}{wx} \right)\Lambda(dx),
\end{equation}
where $\delta_0$ denoted the Dirac measure concentrated in $x=0$. 
It follows that each measure $G_w$ is uniquely determined by the function $\varphi_w$, 
or equivalently
 by the function $\psi$. 
 
 By subtraction, we can suppose $\bar b=0$, $\bar c=0$ and $\bar \Lambda=0$, so that $\psi=0$ and therefore  $\varphi_w=0$ and $G_w(dx)=0$ for every $w \in \R^d$.
 Evaluating the right-hand side of \eqref{Gw} in the singleton  $\{0\}$ we get 
 $$
 w^T c w=0, \quad w \in \R^d. 
 $$
 By the spectral theorem, being $c$ a  symmetric matrix it is (orthogonally) diagonalizable, so  $c= p^T D p$ with $D$ diagonal  and $p$ orthogonal matrix.
 Therefore, setting $\tilde w = p w$, 
 $$
\tilde w^T D \tilde w=0, \quad \tilde w \in \R^d,  
 $$
therefore 
 $D=0$ and $c=0$ as well. 
 
 Going back to \eqref{Gw}, we have that  
 $$
 \left(1- \frac{sin(wx)}{wx} \right)\Lambda(dx)\equiv 0, \quad w \in \R^d. 
 $$
 Since $1- \frac{sin(a)}{a}> 0$ for all $a \neq 0$, and recalling that $\Lambda(\{0\})=0$, we get
 $\Lambda=0$. Finally, from \eqref{245},  we also have $b=0$.
\endproof

 \appendix 
\renewcommand\thesection{Appendix} 
\section{} 
\renewcommand\thesection{\Alph{subsection}} 
\renewcommand\thesubsection{\Alph{subsection}}

\subsection{Jacod-Shiryaev framework}\label{App:A}
For the sake of the reader, we recall below the  equivalence result for c\`adl\`ag semimartingales stated in Theorem  2.42, Chapter II, of  \cite{JacodBook}. 
\begin{theorem}\label{T:equiv_mtgpb_semimart}
Let $X$ be an adapted  c\`adl\`ag process. 	

Let $B^k$ be an  $\check{\mathbb F}$-predictable process, with finite  variation on finite intervals, and $B_0^k=0$, 
 $C$ be
	an $(\check{\mathcal F_t})$-adapted continuous process of finite variation with $C_0=0$, 	and  $\nu$ be an $(\check{\mathcal F_t})$-predictable random measure on $\R_+ \times \R$.
	
There is equivalence between the two following statements.
\begin{itemize} 
\item [(i)] X is a  semimartingale with characteristics $(B^k, C, \nu)$. 
\item [(ii)] For each bounded function $F$ of class $C^{2}_b$,  the process 
\begin{align*}
&F(X_{\cdot}) - F(X_0)  
- \frac{1}{2} \int_0^{\cdot}   F''( X_s) \,d(C\circ X)_s-  \int_0^{\cdot}  F'( X_s) \,d (B^{k}\circ X)\\
 &-  (F(X_{-} + x) -F( X_{-})-k(x)\, F'( X_{-}))\star (\nu\circ X)
\end{align*} 
is a local martingale.
\end{itemize} 
\end{theorem}
{\bf ACKNOWLEDGEMENTS.}
The research of the second named author was partially supported by
the ANR-22-CE40-0015-01 (SDAIM).

\addcontentsline{toc}{chapter}{Bibliography} 
\bibliographystyle{plain} 
\bibliography{../../../../BIBLIO_FILE/BiblioLivreFRPV_TESI} 

\def\cprime{$'$} \def\cprime{$'$} \def\cprime{$'$} \def\cprime{$'$}
\begin{thebibliography}{10}

\bibitem{BandiniRusso1}
E.~Bandini and F.~Russo.
\newblock Weak {D}irichlet processes with jumps.
\newblock {\em {Stochastic Process. Appl.}}, 127(12):4139--4189, 2017.

\bibitem{BandiniRusso2}
E.~Bandini and F.~Russo.
\newblock Special weak {D}irichlet processes and {BSDE}s driven by a random
  measure.
\newblock {\em Bernoulli}, 24(4A):2569--2609, 2018.

\bibitem{BandiniRusso_DistrDrift}
E.~Bandini and F.~Russo.
\newblock Path-dependent {SDE}s with jumps and irregular drift: well-posedness
  and {D}irichlet properties.
\newblock {\em Preprint Arxiv 2211.03444}, 2022.

\bibitem{BandiniRusso_RevisedWeakDir}
E.~Bandini and F.~Russo.
\newblock Weak {D}irichlet processes and generalized martingale problems.
\newblock {\em Stochastic Process. Appl.}, 170(104261), 2024.

\bibitem{BLT}
B.~Bouchard, G.~Loeper, and X.~Tan.
\newblock A {{\(\mathbb{C}^{0, 1}\)}}-functional {It{\^o}}'s formula and its
  applications in mathematical finance.
\newblock {\em Stochastic Process. Appl.}, 148:299--323, 2022.

\bibitem{cjms}
F.~Coquet, A.~Jakubowski, J.~M{\'e}min, and L.~S{\l}omi{\'n}ski.
\newblock Natural decomposition of processes and weak {D}irichlet processes.
\newblock In {\em In memoriam {P}aul-{A}ndr\'e {M}eyer: {S}\'eminaire de
  {P}robabilit\'es {XXXIX}}, volume 1874 of {\em Lecture Notes in Math.}, pages
  81--116. Springer, Berlin, 2006.

\bibitem{nunno-book}
B.~Di~Nunno, G.~{\O}ksendal and F.~Proske.
\newblock {\em Malliavin calculus for {L}\'evy processes with applications to
  finance}.
\newblock Universitext. Springer-Verlag, Berlin, 2009.

\bibitem{er2}
M.~Errami and F.~Russo.
\newblock {$n$}-covariation, generalized {D}irichlet processes and calculus
  with respect to finite cubic variation processes.
\newblock {\em Stochastic Process. Appl.}, 104(2):259--299, 2003.

\bibitem{gr1}
F.~Gozzi and F.~Russo.
\newblock Verification theorems for stochastic optimal control problems via a
  time dependent {F}ukushima-{D}irichlet decomposition.
\newblock {\em Stochastic Process. Appl.}, 116(11):1530--1562, 2006.

\bibitem{gr}
F.~Gozzi and F.~Russo.
\newblock Weak {D}irichlet processes with a stochastic control perspective.
\newblock {\em Stochastic Process. Appl.}, 116(11):1563--1583, 2006.

\bibitem{chineseBook}
S.~He, J.~Wang, and J.~Yan.
\newblock {\em Semimartingale theory and stochastic calculus}.
\newblock Science Press Beijing New York, 1992.

\bibitem{JacodBook}
J.~Jacod.
\newblock {\em Calcul {S}tochastique et {P}robl\`emes de martingales}, volume
  714 of {\em Lecture Notes in Mathematics}.
\newblock Springer, Berlin, 1979.

\bibitem{jacod_book}
J.~Jacod and A.~N. Shiryaev.
\newblock {\em Limit theorems for stochastic processes}, volume 288 of {\em
  Grundlehren der Mathematischen Wissenschaften [Fundamental Principles of
  Mathematical Sciences]}.
\newblock Springer-Verlag, Berlin, second edition, 2003.

\bibitem{rv95}
F.~Russo and P.~Vallois.
\newblock The generalized covariation process and {I}t\^o formula.
\newblock {\em Stochastic Process. Appl.}, 59(1):81--104, 1995.

\bibitem{Russo_Vallois_Book}
F.~Russo and P.~Vallois.
\newblock {\em Stochastic Calculus via Regularizations}, volume~11.
\newblock Springer International Publishing. Springer-Bocconi, 2022.

\bibitem{Stroock_Varadhan}
D.~W. Stroock and S.~S. Varadhan.
\newblock {\em Multidimensional diffusion processes}.
\newblock Springer, 2007.

\end{thebibliography}

\end{document}